\documentclass[letterpaper,10 pt,journal,twoside]{IEEEtran}

\usepackage{amsmath,amsfonts,amssymb,graphicx,graphics,enumerate,dsfont,float,color,bm}

\IEEEoverridecommandlockouts

\newcommand{\argmax}{\mathop{\mathrm{argmax}}}
\newcommand{\maxm}{\mathop{\mathrm{maximize}}}

\newcommand{\Ex}{\mathbb{E}}

\newcommand{\Xstar}{\mathcal{X}^*}
\newcommand{\Astar}{\mathcal{A}^*}

\newcommand{\bbA}{\mathbb{A}}
\newcommand{\bbR}{\mathbb{R}}

\newcommand{\bbPi}{\mathcal{P}}

\def\DelX0{\Delta(X_0)}
\def\delpi0{\delta^{\pi}}

\def\Ex{\textsc{E}}

\newtheorem{thm}{Theorem}

\newtheorem{prop}{Proposition}
\newtheorem{remark}{Remark}

\newcommand{\comment}[1]{}

\title{A General Framework for Bounding Approximate Dynamic Programming Schemes}

\author{Yajing Liu, Edwin K. P. Chong, Ali Pezeshki, and  Zhenliang Zhang
\thanks{This work was supported in part by NSF under awards CCF-1422658, CCF-1712788, and CMMI-1638284, and by the CSU Information Science and Technology Center (ISTeC).}
\thanks{Y. Liu is with National Renewable Energy Laboratory, Golden, CO 80401, USA {\tt\small yajing.liu@nrel.gov}}
\thanks{E. K. P. Chong and A. Pezeshki are with the Department of Electrical and Computer Engineering, and the Department of Mathematics, Colorado State University, Fort Collins, CO 80523, USA {\tt\small Edwin.Chong,Ali.Pezeshki@Colostate.Edu}}
\thanks{Z. Zhang is with Alibaba iDST, Seattle, WA,  {\tt\small zhenliang.zhang@alibaba-inc.com}}
}

\begin{document}
\sloppy

\pagestyle{empty}

\maketitle
\thispagestyle{empty}

\begin{abstract}

For years, there has been interest in approximation methods for solving \emph{dynamic programming} problems, because of the inherent complexity in computing optimal solutions characterized by Bellman's principle of optimality. A wide range of approximate dynamic programming (ADP) methods now exists. It is of great interest to guarantee that the performance of an ADP scheme be at least some known fraction, say $\beta$, of optimal. 
This paper introduces a general approach to bounding the performance of ADP methods, in this  sense, in the stochastic setting. The approach is based on new results for bounding greedy solutions in string optimization problems, where one has to choose a string (ordered set) of actions to maximize an objective function. This bounding technique is inspired by submodularity theory, but submodularity is not required for establishing bounds. Instead, the bounding is based on quantifying certain notions of curvature of string functions; the smaller the curvatures the better the bound. The key insight is that any ADP scheme is a greedy scheme for some surrogate string objective function that coincides in its optimal solution and value with those of the original optimal control problem. The ADP scheme then yields to the bounding technique mentioned above, and the curvatures of the surrogate objective determine the value $\beta$ of the bound. The surrogate objective and its curvatures depend on the specific ADP.

\end{abstract}


\section{Introduction}

\IEEEPARstart{I}{n} sequential decision making, adaptive sensing, and optimal control, we are frequently faced with optimally choosing a string (finite sequence) of actions over a finite horizon to maximize an objective function. In stochastic settings, these problems are often formulated as \emph{stochastic optimal control problems} in the form of Markov decision processes (MDPs) or partially observable Markov decision processes (POMDPs) \cite{Bertsekas12,Vikram_POMDP,ChK09}. A general approach to is to use \emph{dynamic programming} via \emph{Bellman's principle} of optimality (see, e.g., \cite{Bertsekas12,Vikram_POMDP,ChK09}). However, the computational complexity of this approach grows exponentially with the size of the action space and the decision horizon. Because of this inherent complexity, for years there has been interest in developing approximation methods for solving dynamic programming problems, leading to a wide range of approximate dynamic programming (ADP) schemes. These techniques all aim to replace the \emph{expected-value-to-go (EVTG)} term in Bellman's principle, whose computation is intractable, with computationally tractable approximations. Examples of ADP schemes include myopic and rollout policies, reinforcement learning with neural networks, hindsight optimization, foresight optimization, and model-predictive control (see, e.g., \cite{Bertsekas12,ChK09,Ber05,GrR08}). Although a wide range of ADP methods have been developed, in general it is difficult to tell, without doing extensive simulation and testing, if a given ADP scheme has good performance, and even then it is hard to say how far from optimal it is. 

Here, we develop a general framework for bounding the performance of ADP methods relative to the optimal policy in the stochastic setting. By a bound we \emph{specifically} mean a guarantee that the objective value of the ADP scheme is at least a known fraction $\beta$ of the optimal objective value, where $\beta$ depends on the ADP scheme.  A general framework for deriving such bounds for ADP schemes has remained elusive. We note that there are previous results on general performance bounds for ADP schemes, but not of the kind that we seek here. For example, \cite[Props.~3.1 \& 3.2]{Ber05} provides bounds on the difference between the optimal objective value and the one from ADP for the infinite-horizon case, under certain assumptions on the approximation. Another example is \cite{GrR08}, which bounds the difference in the performance between \emph{model-predictive} controllers and the optimal
\emph{infinite-horizon} controller in the deterministic setting. These absolute-difference bounds can be converted into a bound on the absolute difference normalized by the optimal value \cite{GrR08}. However, in general, it is impossible to convert a bound for the normalized difference between two quantities to a bound for their ratio.

Our contribution is different from prior work in several key aspects: (1) We consider finite horizon discrete stochastic optimal control problems; (2) We have a specific notion of bounding that determines \emph{what fraction} of the optimal performance an ADP scheme is guaranteed to achieve; (3) Our bounding approach is based on specific notions of \emph{curvature} for the ADP scheme. 
The practical significance of our contribution is twofold. First, our method provides bounds for an ADP scheme relative to the true optimal performance even though the latter cannot be computed for real-world problems. Second, the bound is
also useful for comparing different ADP schemes without having to do extensive simulation. 

We first derive general lower bounds on the performance of \emph{greedy} solutions to string optimization problems relative to their optimal solutions. By a string optimization problem, we mean a problem in which the objective function is a map from a feasible set of strings (ordered sets) of actions to the real line---the decision variable of the problem is a string. The goal is to select a string of actions to maximize the objective function, subject to a length constraint on the strings (finite decision horizon). This problem is combinatorial; the complexity of finding the optimal string of actions is generally exponential in the cardinality of the action space and the size of the decision horizon. The greedy solution is easy to compute but provides only an approximate solution to the problem. The greedy scheme starts with the empty string and picks at each stage of the optimization an action that maximizes the step-wise (marginal) gain in the objective function. Our bounds, established in Section~\ref{sc:background_bound}, show that any greedy solution is guaranteed to achieve at least a factor $\beta=(1-e^{-\eta(1-\sigma)})/\eta$ of the optimal objective value, where $\eta$ and $\sigma$ denote specific notions of \emph{curvature} of the objective function. The smaller the curvatures the larger the bound. 

The form of the bound discussed above is reminiscent of bounds in submodular optimization problems (see, e.g., 
\cite{nemhauser1978analysis,vondrak2010submodularity}). Our results are indeed \emph{inspired} by our prior work \cite{ZhC16J} on string submodularity---an extension of submodularity theory from functions of sets to functions of strings. However, here we do not need submodularity in deriving our bounds. Whenever submodularity holds the curvature values are generally smaller and the bound becomes larger (better). 

Our main idea is, given an ADP scheme, to formulate a surrogate string optimization problem for the optimal control problem with two properties: (1) The optimal solution and the optimal objective value of the surrogate problem coincide with those of the optimal control problem; (2) The greedy solution to the surrogate problem is the given ADP scheme for the optimal control problem. Then, our framework for bounding greedy solutions of string optimization problems applies to bounding ADP schemes, where the value of the bound depends only on the curvatures of the surrogate objective function. Of course the surrogate objective function and its curvatures depend on the reward function of the optimal control problem and the specific ADP scheme used. In Section~\ref{sc:pr_adp}, we describe how this can be done. 

\section{Stochastic Optimal Control and ADP}\label{sc:background_adp}

\subsection{Stochastic Optimal Control}

In this section, we introduce a general stochastic optimal control problem and explain what approximate dynamic programming (ADP) is.  Our discussion here follows \cite{ChK09}.

We begin with describing a deterministic optimal control problem to establish our notation  and then move to describing a stochastic optimal control problem, which is our focus. Let $\mathcal{X}$ denote a finite set of states and $\mathcal{A}$ a finite set of control actions. Given $x_1\in\mathcal{X}$ and functions $h:\mathcal{X}\times \mathcal{A}\rightarrow \mathcal{X}$ and $g:\mathcal{X}^K\times \mathcal{A}^K\rightarrow\mathbb{R}_{+}$, consider the optimization problem
\begin{align}\label{deterministicproblem}
\begin{array}{l}
\maxm\limits_{a_1,\ldots, a_K\in\mathcal{A}} \ \ g(x_1,\ldots,x_K;a_1,\ldots,a_K) \\
\text{s.\ t.} \ \ x_{k+1}=h(x_k,a_k),\ k=1,\ldots, K-1.
\end{array}
\end{align}
Think of $a_k$ as the \emph{control action} applied at time $k$ and $x_k$ the
\emph{state} visited at time $k$. The function $h$ represents the state-transition law. The real number $g(x_1,\ldots,x_K;a_1,\ldots,a_K)$ is the \emph{total reward} accrued by applying the string of actions $(a_1,\ldots,a_K)$ along the state \emph{path} (trajectory) $(x_1,\ldots,x_K)$. Problem \eqref{deterministicproblem} is called a (deterministic) optimal control problem and the total reward is typically constructed as\vspace{-.1cm}
\begin{equation}\label{eq:gadd}
g(x_1,\ldots,x_K;a_1,\ldots,a_K)=\sum\limits_{k=1}^K r(x_k,a_k),\vspace{-.1cm}
\end{equation}
where $r:\mathcal{X}\times\mathcal{A}\to\mathbb{R}_{+}$ for $k=1,\ldots,K$ is the \emph{immediate reward} accrued at time $k$ by applying action $a_k$ at state $x_k$. From here on, we assume that $g$ has the additive form in \eqref{eq:gadd}. We could have made $g$, $h$, and $r$ explicitly time dependent. However, time can always be incorporated into the state, and so our formulation is without loss of generality.

We now turn our attention to a \emph{stochastic} version of problem~(\ref{deterministicproblem}). The key difference is that the state evolves randomly over time in response to actions. The distribution of states is specified by the state transition law 
$x_{k+1} = h(x_k,a_k,\xi_k)$, $k=1,\ldots,K-1$, where $x_1$ is a given initial state and $\{\xi_k\}_{k=1}^{K-1}$ is an i.i.d.\ random sequence. 
With this modification, we need to change the objective function to 
$\Ex[g(x_1,\ldots,x_K;a_1,\ldots,a_K)|x_1]$, involving \emph{expectation}, where $\Ex[\cdot|x_1]$ represents conditional expectation given the initial state $x_1$. With this specification, the sequence of states $\{x_k\}_{k=1}^K$ has a ``Markovian'' property in the usual sense. This reduces the problem to one of finding, for each time $k$ and each reachable state $x_k^*$, an optimal action $\pi_k^*(x_k^*)$, corresponding to a \emph{state-feedback} control law. This defines a mapping $\pi_k^*$, often called a \emph{policy}  (or, sometimes, a \emph{Markovian policy}). But the chosen action $\pi_k^*(x_k^*)$ is a random variable \emph{adapted to} $\{\xi_{k-1}\}$, meaning it is measurable with respect to $\xi_1,\ldots,\xi_{k-1}$. This type of policy is called a \emph{randomized policy} in the stochastic optimal control literature. Henceforth, we will use the term \emph{policy} to mean randomized policy.

For $k=1,\ldots,K$, let $\pi_k$ be a policy. For convenience, we will also refer to the entire string $(\pi_1,\ldots, \pi_K)$ as simply a policy. The stochastic optimal control problem can be formulated in the following form:
\begin{align}\label{stochasticproblem}
\hspace{-.3cm}\maxm\limits_{\pi_1,\ldots, \pi_K} \ \    \Ex[g(x_1,\ldots,x_K;\pi_1(x_1),\ldots,\pi_K(x_K))|x_1] \nonumber\\
\hspace{-.3cm}\text{s.\ t.} \ \ x_{k+1}=h(x_k,\pi_k(x_k),\xi_k),\ k=1,\ldots, K-1,
\end{align}
where $\Ex[g(x_1,\ldots,x_K;\pi_1(x_1),\ldots,\pi_K(x_K))|x_1]$ equals
\[
\sum\limits_{k=1}^K\Ex[r(x_k,\pi_k(x_k))|x_1],
\]
the conditional expected cumulative reward over a time horizon of length $K$ given the initial state $x_1$. In the stochastic control problem, $(\pi_1,\ldots, \pi_K)$ is the decision variable.

The stochastic optimal control problem \eqref{stochasticproblem} also goes by the name \emph{Markov decision problem (MDP)} (or Markov decision \emph{process}), and arises in a wide variety of areas, including  sensor resource management, congestion control, UAV guidance, and the game of Go (see, e.g., \cite{Bertsekas12,Vikram_POMDP,ChK09,AlphaGo16a}). In problems where the state is only partially observable we will also have an observation law $y_k=c(x_k,\pi_{k-1}(x_{k-1}),\nu_k)$ as a constraint, where $\nu_k$ is the measurement noise at time $k$. In such a setting the state $x_k$ will be replaced with the \emph{belief-state} $b_k$, which is the posterior distribution of the underlying state $x_k$ given the history of observations and actions. The problem is then called a belief-state MDP or a partially observable MDP (POMDP) \cite{Bertsekas12,Vikram_POMDP,ChK09}. From here on we will develop our approach for bounding approximate solutions to MDPs, but all of our results also apply to bounding approximate solutions to POMDPs.

\subsection{Dynamic Programming}

A scheme  or policy $\Pi_K^*=(\pi_1^*,\ldots, \pi_K^*)$ is \emph{optimal} if 
\begin{equation*}
\Pi_K^*\in \argmax_{\pi_1,\ldots,\pi_K} \ \Ex[g(x_1,\ldots,x_K;\pi_1(x_1),\ldots,\pi_K(x_K)|x_1],  
\end{equation*}
where $x_{k+1}=h(x_k,\pi_k(x_k),\xi_k)\ \text{for}\ 1\leq k\leq K-1$ and $\argmax$ is the set of policies that maximize the objective function (there might be multiple possible such optimal policies, hence the notation ``$\in\argmax$''). The optimal policy defined above is characterized by \emph{Bellman's principle} of optimality (also called the dynamic programming principle). To explain, for each $k=1,\ldots,K$, 
let $\bbPi_k$ denote the set of all strings $(\pi_k,\ldots, \pi_K)$. Next, define
functions $V_k: \mathcal{X} \times \bbPi_k \to \mathbb{R}_{+}$ by
$ 
V_k(x_k,\pi_k,\ldots,\pi_K) = \sum_{i=k}^K \Ex[r(x_i,\pi_i(x_i))|x_k],
$
where for $k=1,\ldots, K$, $x_{i+1} = h(x_i,\pi_i(x_i),\xi_i)$, $i=k,\ldots,K-1$.
The objective function of problem \eqref{stochasticproblem} can be written as $V_1(x_1,\pi_1,\ldots,\pi_K)$, where $x_{k+1} = h(x_k,\pi_k(x_k),\xi_k),\ k=1,\ldots,K-1$. Given $x_1$, define $x^*_1=x_1$ and $x^*_{k+1}=h(x^*_k,\pi^*_k(x_k^*),\xi_k)$, $k=1,\ldots,K-1$. Then, Bellman's principle states that 
\begin{equation}\label{eqn:bellman}
\pi_k^*(x^*_k) \in \argmax_{a\in\mathcal{A}} 
Q(x_k^*,a), \ \ k=1,\ldots,K,
\end{equation}
constitute an optimal policy, where  $Q(x_k^*,a)=r(x_k^*,a)+\Ex[V_{k+1}(x^a_{k+1},\pi_{k+1}^*,\ldots,\pi_K^*)|x_k^*,a]$ is  the \emph{$Q$-value} of state $x_k^*$ and action $a$, and  
$x_{k+1}^a=h(x_k^*,a,\xi_k)$ and $x_{i+1}^a=h(x_i^a,\pi_i^*(x_i^a),\xi_i)$ for $i=k+1,\ldots,K-1$,
with the convention that $V_{K+1}(\cdot)\equiv 0$. Moreover, any policy satisfying (\ref{eqn:bellman}) above is optimal. The second term on the right-hand-side of the $Q$-value, $\Ex[V_{k+1}(x^a_{k+1},\pi^*_{k+1},\ldots,\pi^*_K)|x_k^*,a]$, is called the \emph{expected value-to-go} (EVTG). If the problem were a POMDP, the state $x_k^*$ would be replaced by the belief-state $b_k^*$ (the posterior distribution of $x_k^*$ given the past observations and actions). In a POMDP, a policy is a (randomized) mapping on the space $\mathcal{B}$ of belief-states and takes values in $\mathcal{A}$.

Bellman's principle provides a method to compute an optimal solution: We use (\ref{eqn:bellman}) to iterate backwards over the time indices $k=K,K-1,\ldots,1$, keeping the states as variables, working all the way back to $k=1$. This is the familiar \emph{dynamic programming algorithm}. However, the procedure suffers from the \emph{curse of dimensionality} and is therefore impractical for many problems of interest.

\subsection{Approximate Dynamic Programming}\label{sc:background_adpschemes}

An ADP scheme replaces the EVTG term in Bellman’s principle, whose computation is intractable, with a computationally tractable approximation $W_{k+1}(\hat{x}_k,a)$. We start at time $k=1$, at state $\hat{x}_1=x_1$, and for each $k=1,\ldots,K$, we compute the subsequent control actions and states using
\begin{align}
\label{eqn:adp}
\hat{\pi}_k(\hat{x}_k) \in \argmax_{a\in\mathcal{A}} \{r(\hat{x}_k,a) + W_{k+1}(\hat{x}_k,a)\}
\end{align}
and $\hat{x}_{k+1} = h(\hat{x}_k,\hat{\pi}_k(\hat{x}_k),\xi_k)$.
The EVTG approximation term $W_{k+1}(\hat{x}_k,a)$ can be based on a number of methods, including myopic \cite{ChK09}, rollout \cite{BeT97}, reinforcement learning \cite{SuB98}, hindsight/foresight optimization \cite{ChK09}, and model-predictive control \cite{GrR08}. In each of these ADP schemes, the approximation $W_{k+1}$ has a specific form. For example, in the myopic scheme, $W_{k+1}(\hat{x}_k,a)=0$ and the EVTG is simply ignored. In contrast, reinforcement learning uses a parametric function approximator for the EVTG or equivalently the $Q$-value function. The parametric approximator typically is of the form $Q(x,a)=\theta(a)^T \phi(x)$, where $\phi(x)$ is a \emph{feature} vector or \emph{basis function} (often constructed by a domain expert) associated with state $x$ and the coefficients $\theta(a)$ are learned from training data. The usual practice is to use a neural network. Having learned $\theta(a)$, actions are computed according to $\argmax_{a} \theta(a)^\top\phi(x)$. The reader is referred to \cite{ChK09} for expressions of $W_{k+1}$ associated with different ADP schemes.

Our main goal is to develop a general framework for bounding the performance of any ADP scheme relative to the optimal solution. The next
section provides the tool that we will use in Section~IV to develop our bounding framework.

\section{Bounding Greedy Solutions}\label{sc:background_bound}

In this section, we consider string optimization problems, where we wish to maximize an objective function over strings (ordered sets) of actions. We present performance bounds for the greedy solutions to such problems, in terms of certain notions of curvature for the string objective function. Again, by a bound we mean a guarantee that the objective value of the greedy solution is at least a constant factor of the objective value of the optimal solution, where the constant factor is only a function of the curvatures. This discussion is inspired by our prior work on string submodularity \cite{ZhC16J}, where we extended the concept of submodularity, notions of curvature, and associated bounds for greedy solutions, from functions defined on sets to functions defined on strings. But the bounding framework we present here does not require submodularity. 

\noindent\emph{String Optimization:} Let $\bbA$ be a set of possible actions. Let $A=(a_1,a_2,\ldots,a_k)$
denote a \emph{string} of actions taken over $k$ consecutive
stages, where $a_i\in \bbA$ for $i=1,2,\ldots, k$. Let $\bbA^*$
denote the set of all possible strings of actions (of arbitrary
length, including the empty string $\emptyset$) and $f: \bbA^*\to \bbR_{+}$ be an objective function. The goal is to find a string $M\in \bbA^*$, with a length $|M|\le K$ ($K$ prespecified), to maximize $f$:
\begin{align}
\text{maximize } & f(M) \nonumber \\ 
\textrm{subject to } & M\in\bbA^*,\ |M|\leq K.
\label{eqn:2}
\end{align}

\noindent\emph{Monotoneity and Diminishing Return:} Consider two arbitrary strings $M=(a_1,a_2,\ldots, a_m)$ and $N=(b_1,b_2,\ldots, b_n)$ in $\bbA^*$. We define $(M,N)=(a_1,a_2,\ldots, a_m, b_1,b_2,\ldots, b_n)$, as the concatenation of $M$ and $N$. For $M,N\in \bbA^*$, we write $M\preceq N$ if we have $N=(M, L)$ for some $L\in \bbA^*$. In this case, we say that $M$ is a \emph{prefix} of $N$. The function $f$ is said to have the \emph{prefix-monotone} property if for any $M\preceq N \in \bbA^*$ with $|N|\leq K$,  $f(N)\geq f(M)$. Without loss of generality, we assume that $f(\emptyset)=0$. Then $f(M)\geq 0$ holds for any $M\in\bbA^*$ if $f$ is prefix-monotone. Prefix-monotoneity guarantees that the objective function does not decrease by adding a new action. The function $f$ is said to have the \emph{diminishing-return} property if for any $M\preceq N \in \bbA^*$ with $|N|\leq K$ and  $a\in \bbA$, $f((M,a))-f(M) \geq f((N,a))-f(N)$. This property says that the marginal gain of taking any action $a$ early on in the decision horizon is greater than or equal to that of taking the same action later in the decision horizon. It is akin to concavity for functions on the real line.

\noindent\emph{Optimal strategy}: Any solution to \eqref{eqn:2} is called an \emph{optimal} strategy. If $f$ is prefix-monotone, then there exists an optimal strategy with length $K$, denoted $O_K=(o_1,\ldots,o_K)$. 

\noindent\emph{Greedy strategy}: A string $G_K=(b_1,b_2,\ldots,b_{K})$ is called a \emph{greedy} strategy if for all $i=1,2,\ldots,k$,
\begin{align*}
b_i&\in\mathop{\argmax}\limits_{b\in \bbA} f((b_1,b_2,\ldots,b_{i-1},b)).
\end{align*}

\noindent\emph{Curvatures:} In \cite{ZhC16J} and \cite{LiC_ICML18}, we introduced various notions of \emph{curvature}, which measure the extent to which a string function has the diminishing return property, either along particular trajectories in the action space or along all trajectories. The smaller the curvature the greater the extent of diminishing returns. These notions are called curvature for two reasons: one is the analogy between the diminishing return property for string functions and concavity for functions on $\mathbb{R}$; the other is that their expressions are akin to second-order differences. Here we present two specific notions of curvature, which are particularly convenient for our bounding framework. 

For any string $M=(m_1,m_2,\ldots,m_K)\in \bbA^*$, let $M_{i:j}=(m_i,\ldots, m_j)$ for $1\leq i\leq j\leq K$. We define the \emph{total curvature} $\eta$ of $f$ from the greedy trajectory as
\begin{align}\label{eq:etadef}
\eta=&\max\limits_{\substack{M\in\bbA^*\\ 1\leq i \leq K-1}}\frac{K}{K-i}\bigg\{ 1-\nonumber\\
& \quad\quad\quad\quad\quad\frac{f((G_{1:i},M_{i+1:K}))-\frac{K-i}{K}f(M)}{f(G_{1:i})}\bigg\}.
\end{align}

Adding $-f(\emptyset)=0$ to the denominator and completing the fraction reveals that $\eta$ is in fact akin to a second-order difference. Also, if there exists $M_{i+1:K}$ such that $f((G_{1:i},M_{i+1:K}))-f(G_{1:i})\leq (K-i)/K f(O_K)$, then $\eta \geq 0$.

Let $M_{i+1:i}=\emptyset$. We define the \emph{forward curvature} of $f$ from the greedy trajectory as
\begin{align}\label{eq:sigmadef}
\sigma&=\max\limits_{\substack{M_{i+1:j}\in\bbA^*\\0\leq i\leq K-1\\i+1\leq j\leq K}}\
\bigg\{1-\nonumber\\
&\quad\quad\quad\frac{f((G_{1:i},m_j))-f(G_{1:i})}{f((G_{1:i},M_{i+1:j}))-f((G_{1:i},M_{i+1:j-1}))}\bigg\}.
\end{align}

\begin{remark}\label{sigmaremark}
If $f:\mathbb{A}^*\rightarrow\mathbb{R}_{+}$ is prefix-monotone, then $0\leq \sigma\leq 1$. If $f$ has the diminishing-return property, then $\sigma = 0$. The derivations of these results are straightforward and are omitted due to a lack of space.
\end{remark}
 
The following theorem gives a general performance bound for the greedy solution $G_K$, relative to that of the optimal solution $O_K$, in terms of curvatures $\eta$ and $\sigma$. 
\begin{thm}
\label{generalbound1}
If $f: \bbA^*\to \bbR_{+}$ is prefix-monotone, then
\begin{align*}
   \frac{f(G_K)}{f(O_K)}\geq \frac{1}{\eta}\left(1-\left(1-\eta\frac{1-\sigma}{K}\right)^K\right). 
\end{align*}
\end{thm}
\emph{Proof.} See Appendix A.

The bound above is tight \cite{LiC_ICML18} and as $K\to\infty$ converges to $\left(1-e^{-\eta(1-\sigma)}\right)/\eta$ from above.

\noindent\emph{Connection to Submodularity:} The celebrated result of Nemhauser \emph{et al.}~\cite{nemhauser1978analysis} states that for maximizing a monotone submodular function over a uniform matroid, the objective value of the greedy strategy is no less than a factor $(1-e^{-1})$ of that of the optimal strategy. Sharper bounds of the form $(1-e^{-\gamma})/\gamma$, with $0\leq \gamma\leq 1$, involving a notion of curvature $\gamma$ for set submodular functions were established in \cite{vondrak2010submodularity}. The concept of submodularity was extended to functions defined over strings in \cite{ZhC16J,golovin2011adaptive}, leading to similar bounds (with and without curvature) in sequential optimization problems. For a survey of bounds involving submodularity, see \cite{LiC20DEDS}.

The bounds in Theorem~\ref{generalbound1} are similar to those from submodularity theory. But submodularity is in fact not needed for deriving them, as we have shown in Theorem~\ref{generalbound1}. Such bounds, in terms of \emph{properly defined} notions of curvature and subject to a monotoneity condition, always hold. When $\eta=1$ and $\sigma=0$, the second bound in Theorem~\ref{generalbound1} is $(1-e^{-1})\approx 0.63$, coinciding with the Nemhauser bound. When  submodularity holds ($\eta<1$ and $\sigma=0$) the bound is better than $(1-e^{-1})$, with practical impact \cite{SuC19}. A difference between our notions of curvature and other notions used in submodularity literature is that the values of our curvatures, and therefore the bounds, do not depend on the behavior of the objective function on strings longer than $K$ (the decision horizon), whereas other curvatures do depend on the values of the set/string function on larger sets/strings. This is a subtle but important difference, because performance bounds for an optimization problem over a horizon of size $K$ should not depend on what the objective function does beyond the decision horizon.

\section{Bounding ADP Schemes}\label{sc:pr_adp}

Our idea is to formulate a surrogate string optimization problem over the set of policy strings such that (1) its greedy solution coincides with the ADP scheme, and (2) its optimal solution and optimal objective value coincide with those of the stochastic optimal control problem. This enables us to bound the performance of the ADP scheme relative to the optimal scheme, in terms of the curvatures of the surrogate objective function. This surrogate function of course must depend on both the immediate reward function $r$ and the approximation $W$ to the EVTG. The key to establishing this result is a fundamental connection between two classes of approximate solutions to general stochastic control problems.
\subsection{PDAO Schemes and GPS Schemes}\label{sc:pr_pdao}

Let $\Xstar=\mathcal{X}\cup\mathcal{X}^2\cup\cdots$ denote the collection of all strings of states and  $\Astar=\mathcal{A}\cup\mathcal{A}^2\cup\cdots$ the collection of all strings of actions. Let $\tilde{g}:\Xstar\times\Astar\rightarrow\mathbb{R}_{+}$ be an objective function. Consider the stochastic control problem 
\begin{align}\label{eq:stochasticproblemf}
\hspace{-.17cm}\maxm\limits_{\pi_1,\ldots, \pi_K} \ \    \Ex[\tilde{g}(x_1,\ldots,x_K;\pi_1(x_1),\ldots,\pi_K(x_K))|x_1] \nonumber\\
\hspace{-.6cm}\text{s.\ t.} \ \ x_{k+1}=h(x_k,\pi_k(x_k),\xi_k),\ k=1,\ldots, K-1.
\end{align}
We are distinguishing between the objective function $g$ of \eqref{stochasticproblem}, which is a function of $K$ states and $K$ actions, and the function $\tilde{g}$ above, which can take arguments with state and action strings that are of arbitrary length. Later on, we connect $\tilde{g}$ to the optimal control objective $g$ and the ADP approximation $W$ in a very specific way to establish our bounding technique for ADPs. Below, we introduce two classes of approximate solutions to \eqref{eq:stochasticproblemf}. We assume throughout that $x_1\in\mathcal{X}$ is given.

Given $x_1^{p}=x_1$, the policy $\Pi_K^p=(\pi_1^{p},\ldots,\pi_K^{p})$ is called a \emph{path-dependent action optimization (PDAO)} scheme if for $k=1,\ldots, K$,
\begin{align}
&\pi_k^{p}(x_k^{p})\in\argmax_{a\in\mathcal{A}} \tilde{g}(x_1^p,\ldots,x_{k}^{p};\pi_1^p(x_1^p),\ldots,\pi_{k-1}^p(x_{k-1}^p),a),\label{greedy1}
\end{align}
where $x_{i+1}^{p}=h(x_i^{p},\pi_i^{p}(x_i^{p}),\xi_i)$ for $1\leq i\leq k-1$.

Given $x_1^{g}=x_1$, the policy $\Pi_K^g=(\pi_1^{g},\ldots,\pi_K^{g})$ is called a \emph{greedy policy-selection (GPS)} scheme if for $k=1,\ldots, K$, 
\begin{align}
\pi_k^{g}\in
\argmax_{\pi_k}\, &\Ex[\tilde{g}(x_1^g,\ldots,x_{k}^{g};\nonumber\\
&\pi_1^g(x_1^g),\ldots,\pi_{k-1}^g(x_{k-1}^g),\pi_k(x_{k}^{g}))|x_1],\label{greedy2}
\end{align}
where $x_{i+1}^{g}=h(x_i^{g},\pi_i^{g}(x_i^{g}),\xi_i)$ for $1\leq i\leq k-1$.

At each time $k$, a PDAO scheme chooses an action based on the \emph{sample path} $\xi_1,\ldots,\xi_{k-1}$, and the chosen action is adapted to $\{\xi_{k-1}\}$. In contrast, a GPS scheme chooses a policy at each time $k$ based on the expected reward. Nonetheless, a PDAO scheme still defines a particular policy. A key result for our approach is the following.

\begin{thm}
\label{theorem3}
Any PDAO scheme is also a GPS scheme: Given a PDAO policy $\Pi_K^p=(\pi_1^{p},\ldots, \pi_K^{p})$, satisfying (\ref{greedy1}), there exists a GPS policy $\Pi_K^g=(\pi_1^{g},\ldots, \pi_K^{g})$ such that $\pi_j^{p}=\pi_j^{g}$ for $1\leq j\leq k$.
\end{thm}

\emph{Proof:} See Appendix B.

\subsection{Bounding PDAO Schemes}

Let $\bbPi^*$ be the set of all strings of policies $(\pi_1,\ldots,\pi_k)$ with $k=0,1,2,\ldots\,$; the case $k=0$ corresponds to the empty string.
Given $x_1$, define 
$\tilde{g}_{\text{avg}}:\bbPi^*\rightarrow\mathbb{R_+}$ by
\begin{align*}
\tilde{g}_{\text{avg}}(\pi_1,\ldots,\pi_k)=\Ex[\tilde{g}(x_1,\ldots,x_k;\pi_1(x_1),\ldots,\pi_k(x_k))|x_1].
\end{align*} 
It is clear that $\tilde{g}_{\text{avg}}(\pi_1,\ldots,\pi_K)$, with $k=K$, is the objective function in \eqref{eq:stochasticproblemf}. So we can convert \eqref{eq:stochasticproblemf} to the following optimization problem, where the objective function is simply a function of policy strings:
\begin{align}
\text{maximize } & \tilde{g}_{\text{avg}}(M) \nonumber \\ 
\textrm{subject to } & M\in\bbPi^*,\ |M|\leq K.
\label{eqn:s2}
\end{align}
Naturally, optimal solutions for the two problems are identical. Moreover, the GPS scheme for \eqref{eq:stochasticproblemf} coincides with the greedy solution to \eqref{eqn:s2}. Therefore, provided $\tilde{g}_{\mathrm{avg}}$ is prefix-monotone, the GPS scheme $\Pi_K^g$ can be bounded as in Theorem~\ref{generalbound1} in terms of the curvatures $\eta$ and $\sigma$ of $\tilde{g}_{\text{avg}}$. At the same time, Theorem~\ref{theorem3} established that any PDAO scheme is a GPS scheme. Hence, we have the following. 
\begin{thm}
\label{thm4}
Let $\tilde{\Pi}_K^*=(\tilde{\pi}_1^*,\ldots, \tilde{\pi}_K^*)$ denote an optimal solution to $\eqref{eq:stochasticproblemf}$.
If $\tilde{g}_{\mathrm{avg}}$ is prefix-monotone, then 
any PDAO scheme $\Pi_K^p=(\pi_1^p,
\ldots, \pi_K^p)$ for problem~(\ref{eq:stochasticproblemf}) achieves the following bound:
\begin{align*}
\frac{\tilde{g}_{\mathrm{avg}}(\Pi_K^p)}{\tilde{g}_{\mathrm{avg}}(\tilde{\Pi}_K^*)}&\geq \frac{1}{\eta}\left(1-\left(1-\eta\frac{1-\sigma}{K}\right)^K\right)>\frac{1-e^{-\eta(1-\sigma)}}{\eta},
\end{align*}
where $\eta$ and $\sigma$ are curvatures of $\tilde{g}_{\mathrm{avg}}: \bbPi^*\rightarrow \mathbb{R}$ as defined in \eqref{eq:etadef} and \eqref{eq:sigmadef}, respectively.
\end{thm}

\subsection{Bounding ADP Schemes}\label{sc:pr_badp}

What does bounding PDAO schemes have to do with bounding ADP schemes? We show here that any ADP scheme is the PDAO scheme for a surrogate function $\tilde{g}$, whose expected value $\tilde{g}_{\mathrm{avg}}$, given $x_1$, is equal to the objective function of our stochastic optimal control problem in \eqref{stochasticproblem}. Therefore, Theorem~\ref{thm4} applies to bounding the ADP scheme. Indeed, define the function $\tilde{g}:\Xstar\times\Astar\to\mathbb{R}_+$ by,\vspace{-.1cm} 
\begin{equation}\label{eq:gtilde}
\tilde{g}(x_{1:k},\Pi_{k}(x_{1:k})) =  \sum_{i=1}^k r(x_i,\pi_i(x_i)) + W_{k+1}(x_k,\pi_k(x_k))\vspace{-.1cm}
\end{equation}
for $k=1,\ldots,K$,
where $x_{1:k}=(x_1,\ldots,x_k)$, $\Pi_{k}(x_{1:k})=(\pi_1(x_1),\ldots,\pi_k(x_k))$,
$x_{k+1}=h_k(x_k,\pi_k(x_k),\xi_k)$, and $W_{K+1}(\cdot)\equiv 0$
by convention. For this $\tilde{g}$, we have an associated PDAO scheme. At the terminal stage $k=K$, by the definition of $\tilde{g}_{\text{avg}}$, we have $\tilde{g}_{\text{avg}}((\pi_1,\ldots,\pi_K))
=\Ex[\tilde{g}(x_{1:K},\Pi_{K}(x_{1:K}))|x_1]=\sum_{i=1}^K\Ex[ r(x_i,\pi_i(x_i))|x_1]$. 
This is equal to the objective function for the original stochastic optimal problem \eqref{stochasticproblem}, and is also the function to be maximized at the final stage of the GPS scheme. By Theorem~\ref{theorem3}, the PDAO scheme associated with the above surrogate $\tilde{g}$ is the GPS scheme for the optimal control problem \eqref{stochasticproblem}. Next, notice that the PDAO scheme, denoted here by $(\hat{\pi}_1(\hat{x}_1),\ldots,\hat{\pi}_{k-1}(\hat{x}_{k-1})$, is given by
\begin{align*}
\hat{\pi}_k(\hat{x}_k)\in \argmax_{a\in\mathcal{A}} \ &  \tilde{g}(\hat{x}_1,\ldots,\hat{x}_{k-1},\hat{x}_k;\\
&\quad\quad\hat{\pi}_1(\hat{x}_1),\ldots,\hat{\pi}_{k-1}(\hat{x}_{k-1}),a)
\nonumber.
\end{align*}
Substituting for $\tilde{g}$ from \eqref{eq:gtilde} and dropping the term $\sum_{i=1}^{k-1} r(\hat{x}_i,\hat{\pi}_i(\hat{x}_i))$, which does not depend on $a$, we have that $
\hat{\pi}_k(\hat{x}_k)= \argmax_{a\in\mathcal{A}} \{r(\hat{x}_k,a) + W_{k+1}(\hat{x}_k,a)\}$. But this is simply the ADP scheme in (\ref{eqn:adp}).

\begin{prop}
The ADP scheme in (\ref{eqn:adp}) is a PDAO scheme for the optimization problem defined above.
\end{prop}

Finally, it remains to establish a sufficient condition for $\tilde{g}_{\text{avg}}$ to be prefix-monotone. A simple calculation using \eqref{eq:gtilde} shows that $\tilde{g}_{\mathrm{avg}}$ is prefix-monotone if, for any $(\pi_1,\ldots,\pi_m)\preceq(\pi_1,\ldots,\pi_n)\in\bbPi^*$ with $1\leq m\leq n\leq K$, we have
$
\Ex[W_{m+1}(x_m, \pi_m(x_m))-W_{n+1}(x_n, \pi_n(x_n))] \leq \sum_{i=m+1}^n\Ex[ r(x_i,\pi_i(x_i))|x_1]$.
Monotoneity holds if the marginal reduction in the approximation to the EVTG over any $n-m$ consecutive steps is no greater than the cumulative expected reward over those steps. This condition is trivially satisfied for the myopic heuristic ($W_{m+1}=0$ for all $m$), if all immediate rewards are chosen to be non-negative (which can always be achieved). For other ADP schemes, the condition sets a constraint on the relation between values of immediate reward and those of the approximation to the EVTG.

Having established a condition for monotoneity and that any ADP scheme is a PDAO scheme, Theorem~\ref{thm4} can be used to bound the performance of any ADP scheme. This bounding framework also guides the design of good ADP schemes, namely by designing the approximate EVTG term $W_{k+1}$ such that the corresponding $\tilde{g}_{\text{avg}}$ has small curvatures. This is a subject for our immediate future work, along with establishing ways to compute or bound the curvatures $\eta$ and $\sigma$ with polynomial number of function evaluations.

\section{Appendix A: Proof of Theorem 1}

For simplicity let $G_i=G_{1:i}$, $O_i=O_{1:i}$. First, we prove 
\begin{equation}
\label{importineq1}
f(G_1)\geq \frac{1-\sigma}{K}f(O_K).
\end{equation}
By the definition of $\sigma$ and the prefix-monotone property of $f$, we have, for $j=2, \ldots, K$, $f((o_j))\geq (1-\sigma)(f(O_j)-f(O_{j-1}))$. This is obtained from $\eqref{eq:sigmadef}$ by observing that $\sigma$ is a max over the set of all strings $M_{i+1:j}$, $1\leq i+1\leq j\leq K$. So $M_{i+1:j}=O_j$
with $i=0$, is an element of that set. 
Summing over $j$ gives 
$\sum_{j=2}^Kf((o_j))+(1-\sigma)f((o_1))\geq (1-\sigma)f(O_K)$. By Remark~\ref{sigmaremark}, we have that $\sigma\geq 0$,
which implies that $1-\sigma\leq 1$.
By definition of the greedy strategy, we have that $f(G_1)\geq f((o_j)) \ \text{for}\ 2\leq j\leq K$. Combining this with the previous inequality and $1-\sigma\leq 1$ gives (\ref{importineq1}).

Second, we prove that for $1\leq i\leq K-1$,
\begin{equation}
\label{importineq2}
f(G_{i+1})\geq \frac{1-\sigma}{K}f(O_K)+\left(1-\eta\frac{1-\sigma}{K}\right)f(G_i).
\end{equation}
Again from the definition of $\sigma$ (with a similar max argument), due to the prefix-monotone property of $f$, for a fixed $i$ with $2\leq i+1\leq j\leq K$, we have that $f((G_i,m_j))-f(G_i)\geq(1-\sigma)(f((G_i,M_{i+1:j}))-f((G_i,M_{i+1:j-1})))$. Summing this over $j$ from $i+1$ to $K$ results in 
$\sum_{j=i+1}^K (f((G_i,m_j))-f(G_i))\geq(1-\sigma)(f((G_i,M_{i+1:K}))-f(G_i))$. For the greedy strategy, we have $f(G_{i+1})-f(G_i)\geq f((G_i,m_j))-f(G_i)$. Combining this with the previous inequality yields 
\begin{equation}
\label{deriveimp1}
f(G_{i+1})-f(G_i)\geq\frac{1-\sigma}{K-i}(f((G_i,M_{i+1:K})))-f(G_i)).
\end{equation}

Using the definition of the curvature $\eta$ in \eqref{eq:etadef}, letting $M=O_K$ (and hence $M_{i+1:K}=O_{i+1:K})$ and noting the max operation, we obtain
\begin{equation}
\label{deriveimp2}
f((G_i,O_{i+1:K}))-f(G_i)\geq \frac{K-i}{K}(f(O_K)-\eta f(G_i)).
\end{equation}
Combining (\ref{deriveimp1}) and (\ref{deriveimp2}), observing that (\ref{deriveimp1}) holds for any $M_{i+1:K}\in\mathbb{A}^\ast$ (including $M_{i+1:K}=O_{i+1:K})$, yields (\ref{importineq2}).

Applying \eqref{importineq2} successively from $i=K-1$ to $i=1$ gives 
\begin{align}\label{importineq3}
f(G_K)&\geq \frac{1-\sigma}{K}f(O_K)+\left(1-\eta\frac{1-\sigma}{K}\right)\frac{1-\sigma}{K}f(O_K)\nonumber\\
 &\quad\quad+\cdots+\left(1-\eta\frac{1-\sigma}{K}\right)^{K-1}f(G_1).
\end{align}
Applying \eqref{importineq1} to the right hand side of \eqref{importineq3} yields
\[
f(G_K)> \frac{1}{\eta}\left(1-\left(1-\eta\frac{1-\sigma}{K}\right)^{K} \right)f(O_K).
\]

\section{Appendix B: Proof of Theorem 2}

Suppose that we are given a PDAO policy $(\pi_1^{p},\ldots, \pi_K^{p})$ (i.e., satisfying (\ref{greedy1})). We will show that there exists a GPS policy $(\pi_1^{g},\ldots, \pi_K^{g})$ such that the two policies are equal, i.e., $\pi_j^{p}=\pi_j^{g}$ for $1\leq j\leq k$. We show this by induction on $k$. For $k=1$, by (\ref{greedy1}), we have that for any $\pi_1$, $\tilde{g}(x_1^{p},\pi_1^{p}(x_1^{p}))\geq \tilde{g}(x_1^{p},\pi_1(x_1^{p}))$, which implies that $\Ex[\tilde{g}(x_1^{p},\pi_1^{p}(x_1^{p}))|x_1]\geq\Ex[\tilde{g}(x_1^{p},\pi_1(x_1^{p}))|x_1]$. Because $x_1^{p}=x_1$, this shows that $\pi_1^p=\pi_1^g$. For the induction step, assume that there exists $(\pi_1^{g},\ldots, \pi_k^{g})$ satisfying (\ref{greedy2}) such that  $\pi_j^{p}=\pi_j^{g}$ for $1\leq j\leq k$. To complete the proof, it suffices to show that $\pi_{k+1}^p$ satisfies (\ref{greedy2}). By definition, $x_{j+1}^{p}=h_j(x_j^{p},\pi_j^{p}(x_j^{p}),\xi_j)$ and 
$x_{j+1}^{g}=h_j(x_j^{g},\pi_j^{g}(x_j^{g}),\xi_j)$ for $1\leq j\leq k$. By the assumption that $\pi_j^{p}=\pi_j^{g}$ for $1\leq j\leq k$ and $x_1^{p}=x_1^{g}$, we have that $x_{j+1}^{p}=x_{j+1}^{g}$ for $1\leq j\leq k$. Thus,  $x_{k+1}^{p}=x_{k+1}^{g}$. For $\pi_{k+1}^{p}$, by (\ref{greedy1}), we have that for any $\pi_{k+1}$, $\tilde{g}(x_1^p,\ldots,x_{k+1}^{p},\pi_1^p(x_1^p),\ldots,\pi_{k+1}^{p}(x_{k+1}^{p}))
\geq 
\tilde{g}(x_1^p,\ldots,x_{k+1}^{p},\pi_1^p(x_1^p),\ldots,\pi_{k+1}(x_{k+1}^{p}))$. This implies that
\begin{align*}
\Ex[(\tilde{g}(x_1^p,\ldots,x_{k+1}^{p},\pi_1^p(x_1^p),\ldots,\pi_{k+1}^{p}(x_{k+1}^{p}))|x_1]\ge\quad\quad\quad\\ \quad\quad\quad\Ex[\tilde{g}(x_1^p,\ldots,x_{k+1}^{p},\pi_1^p(x_1^p),\ldots,\pi_{k+1}(x_{k+1}^{p}))|x_1]. 
\end{align*}
Because $x_{k+1}^{p}=x_{k+1}^{g}$, this means that $\pi_{k+1}^{p}$ satisfies (\ref{greedy2}).
This completes our induction argument. \hfill 

\bibliography{reference}
\bibliographystyle{IEEEtran}

\end{document}